\documentclass[reqno]{amsart}

\usepackage{graphicx,amsmath,amsfonts}
\usepackage{amssymb,amscd,color}
\usepackage[mathscr]{eucal}

\theoremstyle{plain}
\newtheorem{prob}{Problem}

\newtheorem{theorem}{Theorem}[section]
\newtheorem{lemma}[theorem]{Lemma}

\theoremstyle{definition}
\newtheorem{definition}[theorem]{Definition}

\theoremstyle{remark}
\newtheorem*{remark}{Remark}
\newtheorem*{algo}{Algorithm}
\newtheorem*{note}{Note}

\newcommand{\CC}{\hbox{\ensuremath{\mathcal{C}}}}

\newcommand{\N}{\hbox{\ensuremath{\mathbb{N}}}}
\newcommand{\Z}{\hbox{\ensuremath{\mathbb{Z}}}}

\newcommand{\R}{\hbox{\ensuremath{\mathbb{R}}}}

\newcommand{\BB}{\hbox{\ensuremath{\mathfrak{B}}}}
\newcommand{\DD}{\hbox{\ensuremath{\mathfrak{L}}}}

\newcommand{\proj}{\mathscr P}

\newcommand{\HH}{\mathcal H}
\newcommand{\LL}{\mathcal L}
\newcommand{\SiS}{\mathcal S}
\newcommand{\VV}{\mathfrak{B}}
\newcommand{\vv}{\mathcal W}

\newcommand{\EV}{{\bf V}}
\newcommand{\ES}{{\bf S}}
\newcommand{\F}{\hbox{\ensuremath{\mathcal{F}}}}

\newcommand{\zo}{\ensuremath{\omega}}

\newcommand{\sspan}{\text{span}}

\newcommand{\Pil}{\Omega_l}

\newcommand{\beq}{\begin {equation}}
\newcommand{\eeq}{\end{equation}}

\newcommand{\SIS}{shift-invariant space}

\newcommand{ \word}{bundle}

\DeclareMathOperator{\CW}{G}
\DeclareMathOperator{\CP}{H}

\numberwithin{equation}{section}
\begin{document}

\title{Optimal Non-Linear  Models for Sparsity and Sampling}

\author[Akram Aldroubi]{Akram~Aldroubi} \address{Department of Mathematics\\
Vanderbilt University\\ 1326 Stevenson Center\\ Nashville, TN 37240}
\email[Akram Aldroubi]{akram.aldroubi@vanderbilt.edu}

\author[Carlos Cabrelli]{Carlos~Cabrelli}
\address[C. Cabrelli and U. Molter]{Departamento de
Matem\'atica \\ Facultad de Ciencias Exactas y Naturales\\ Universidad
de Buenos Aires\\ Ciudad Universitaria, Pabell\'on I\\ 1428 Capital
Federal\\ ARGENTINA\\ and CONICET, Argentina}
\email[Carlos~Cabrelli]{cabrelli@dm.uba.ar} \thanks{The research of
Akram Aldroubi is supported in part by NSF Grant DMS-0504788. The research of
       Carlos Cabrelli and Ursula Molter is partially supported by
Grants: PICT 15033, CONICET, PIP 5650, UBACyT X058 and X108}

\author[U.Molter]{Ursula~Molter}
\email[Ursula~M.~Molter]{umolter@dm.uba.ar}

\keywords{Sampling, Sparsity, Compressed Sensing, Frames}
\subjclass[2000]{Primary 41A65, 42C15 ; Secondary 68P30, 94A20}

\date{\today}

\begin{abstract}
Given a set of vectors (the data) in a Hilbert space $\HH$, we prove the existence of an optimal collection of  subspaces minimizing the sum of the square of the distances between each vector and its closest subspace in the collection. 
This collection of subspaces gives the best sparse representation for the given data,
in a sense defined in the paper,
and provides an optimal model for sampling in union of subspaces.
The results are proved in a general setting and then applied to the case of low dimensional subspaces of $\R^N$ and to  infinite dimensional shift-invariant spaces in $L^2(\R^d)$. 
We also present an iterative search algorithm for finding the solution subspaces.
These results are tightly  connected to the new emergent theories of compressed sensing and dictionary design,  signal models for signals with finite rate of innovation, and the subspace segmentation problem.
\end{abstract}

\maketitle


\section{introduction}

A new paradigm for signal sampling and reconstruction recently developed by Lu and Do \cite {LD07} starts from the point of view that signals live in some
union of subspaces $\mathcal{M}=\cup_{i\in I}V_i$, instead of a single vector
space $\mathcal{M}=V$ such as the space of band-limited functions also known as the Paley-Wiener space.  This new paradigm is general and includes (when $\mathcal{M}=V$) the classical Shannon sampling theory and its extensions \cite {AG01}, as well as sampling of signal with finite rate of innovation  (see e.g., \cite {MV05,DVB07}). In the
new framework, when we have more than one subspace, the signal space model $\mathcal{M}=\cup_{i\in I}V_i$ is non-linear
and the techniques for reconstructing a signal  $f\in \cup_{i\in I}V_i$ from its samples $\{f(x_j)\}_{j}$ are involved
and the reconstruction operators are non-linear.  

Since for each class of signals  the starting point
of this new theory is the knowledge of the signal space $\mathcal{M}=\cup_{i\in
I}V_i$, the first step for implementing the theory is to find 
an appropriate signal model $\mathcal{M}=\cup_{i\in I}V_i$ from a set of
observed data $\F = \{f_1, \dots, f_m\} $. For
the classical sampling theory, the problem of finding the shift-invariant space model $\mathcal{M}=V$ from a set of observed data has been
studied and solved in \cite {ACHMR03},\cite {ACHM07}. For the new sampling paradigm, the problem consists  in proving the existence and finding subspaces
$V_1,\cdots,V_l,$  {  of some Hilbert space $\HH$} that minimize
the expression
 \begin{equation}
 \label{total-error-0}
e(\F,\big \{V_1,\dots,V_l \big \}) = \sum_{i=1}^m \min_{1\leq j \leq l} d^2(f_i,V_j),
\end{equation}
{  over all possible choices of $l$ subspaces belonging to an appropriate class of subspaces of $\HH$.} Here $\F = \{f_1, \dots, f_m\} \subset \HH$ is a set of observed data and
$d$ is the distance function in $\HH$.

It is well known that the problem of sampling and reconstruction of signals with finite rate of innovation is closely related to the developing theory of compressed sensing (see e.g.,   \cite {CRT06,CR06,CT06,Dev07,Don06,RSV06} and the references therein). Compressed sensing proposes to find a vector $x\in \R^N$ from the knowledge of the values, when applied to $x$, of a relatively small set of functionals $\{\psi_k: \; k=1,\dots,p\}$ (where $p<<N$).  Obviously, the problem of finding $x$ from the set $\{y_k=\langle x, \psi_k\rangle: \, k = 1,\dots,p\}$ is ill-posed. However, it becomes meaningful if $x$ is assumed to be sufficiently sparse. 

A typical assumption of sparsity is that $x$ has at most $n$ non-zero components ($\|x\|_0\le n$), where $n\le 2p<<N$. As a consequence of this assumption of sparsity,  the vector $x$ belongs to some union of subspaces, each of which is generated by exactly $n$ vectors from the canonical basis of $\R^N$. In matrix formulation this problem can be stated as follows: find $x\in \R^N$ with $\|x\|_0\le n$ from the matrix equation $y=Ax$ where $A$ is a $p\times N$ matrix and $y$ is a given vector in $\R^p$. 

A related problem consists in finding an approximation to the vector $y$ using a sparse vector $x$. Formally, this problem can be stated as follows:  find $\min\limits_{x} \|x\|_0$ subject to the constraint $\|Ax-y\|_2\le \varepsilon$ for some given $\varepsilon$. The above two problems, their analysis, extensions, and efficient algorithms for finding their solutions can be found in \cite {AEB06b,AEB06a,BDDW07,CRT06,CR06,CT06,Dev07,Don06,GN03,Tro04} and the references therein.  

If in the above problems the matrix $A$ is also an unknown to be found together with the set of unknown vectors $\{x_i: \; i=1,\dots,m\}\subset \R^N$, then these problems become the problems of finding a {\it dictionary} $A$ from the data $\{y_i: \; i=1,\dots,m \} \subset \R^p$ obtained by sampling the sparse vectors $\{x_i: \; i=1,\dots,m\}\subset \R^N$ see e.g., \cite {AEB06a,AEB06b,GN03}.   In this context, the columns of $A$ are called {\it atoms}  of $A$. Under appropriate assumptions on the data and dictionary, the problem has a unique solution up to a permutation of the columns of $A$ \cite {AEB06a, AEB06b}. Finding the solution to this problem by exhaustive methods is computationally intractable,  but the K-SVD algorithm described in  \cite {AEB06b} provides a computationally effective search algorithm. 

The problem  of finding the signal model for signals with  finite rate of innovation consists of finding a set $\mathcal{M}=\cup_{i\in I}V_i$, formed by subspaces $V_i$ that are infinite dimensional, in general, but  usually structured, e.g., each $V_i$ is a shift-invariant space.  However,  the signal modeling problem  as described by \eqref {total-error-0} is closely related to the dictionary design problem for sparse data, described in the previous paragraph.  

To see this relation, let us formulate the dictionary design  problem as follows:
 given  a class of signals, determine if there exists a dictionary of small size, such that each of the signals can be represented with minimal sparsity.
 
More precisely, assume that we have a class of $m$ signals, where $m$ is a very large number.
We want to know whether there exists a dictionary, such that every signal in the class
is a linear combination of at most $n$ atoms in the dictionary.
Clearly, to make the problem meaningful and realistic the length of the dictionary should be small compared with $m$.

It follows, that if for a given set of data such a dictionary exists,  then the data can be partitioned
into subsets each of which  belongs to a subspace of dimension at most $n$ (i.e. to the subspace generated by the atoms that the signal uses in its representation). That is, each subset of the partition can be associated  to a low dimensional subspace.

Conversely, if our class of signals can be partitioned into $l$ subsets, such that the signals in each subset
belong to a subspace of dimension no bigger than $n$, then by choosing a set of generators from each of the subspaces, we can construct a dictionary of length at most $ln$ with the property that each of the signals can be represented using at most $n$ atoms in the dictionary.

This suggests that the problem of finding a dictionary where the signals have sparse representation
can be solved by finding a small collection of low dimensional subspaces containing our signals,
and viceversa.

So, we will say that the class of signals is {\it $(l,n)$-sparse} (see also Definition~\ref{def-sparse}) if there exist $l$ subspaces
 of dimension at most $n,$ such that the signals in our class belong to the union of these $l$ subspaces. From the above discussion, it is clear that if our data is $(l,n)$-sparse then there exists a dictionary of length at most $ln$. 

A related problem is the subspace segmentation problem for a set of signals in $\R^N$ (see for example \cite{MDHW07,MYDF07}).  This problem occurs in the context of segmentation clustering and classification, and consists in finding whether there exist $l$ subspaces of dimension at most $n$, such that the signals in the class belong to the union of these $l$ subspaces. The subspace segmentation problem has important applications in computer vision, image processing and other areas of  engineering, and it has recently been solved using algebraic methods and algebraic geometrical tools \cite {VMS05}. The method for solving it (known as the Generalized Principle Component Analysis (GPCA)) has also been extended to deal with moderate noise in the data \cite {VMS05}. Moreover, the uniqueness problem has been addressed in \cite{MYDF07}.

Now assume that for a given $l$ and $n$ our data is not $(l,n)$-sparse. {  In that case we prove that there still exists a collection of optimal subspaces providing the needed sparsity. More precisely, if $\varepsilon > 0$ is given, we determine that  there exists a collection} of $l$ subspaces  of dimension at most $n$ such that  the sum of the squares  of the distance of each signal to the union of the subspaces (i.e., the {\em total error}) is not larger than $\varepsilon$, (see formula \eqref{total-error-0}).
 In that case we will say that our data is {\it $(l,n,\varepsilon)$-sparse}. 
 
 As before it is clear that if our data is $(l,n,\varepsilon)$-sparse,  then
 a dictionary of length at most $ln$ exists such that every signal in our class can be approximated using a linear combination of at most $n$ atoms from the dictionary, with total error not larger than $\varepsilon$.
 
Note that this definition of sparsity is an intrinsic property of the data and the space where they belong to, and does not depend on any fixed dictionary.

A relevant and important question is then, given a class of signals and a small number $n$,
which is the minimun possible $\varepsilon$ such the data is $(l,n,\varepsilon)$-sparse? 

In this paper we present a general scheme that allows us to solve the problem described in \eqref {total-error-0}, thereby finding the signal model for the new signal sampling paradigm described in \cite {LD07}, finding a new method for solving the segmentation subspace problem that is optimal in the presence of noise \cite {VMS05},   and solving the $(l,n,\varepsilon)$-sparsity problem (in the sense defined above) for a given set of data, in different contexts. Specifically, given a set $\F$ of $m$  vectors and  numbers $l,n$ such that $ n, l <m $, we prove the existence of no more than $l$  subspaces of dimension no bigger than $n$ that provide  the minimum $\varepsilon$ such that the  vectors in $\F$ are $(l,n,\varepsilon)$-sparse. When the minimum $\varepsilon$ is zero, the data is $(l,n)$-sparse. We also give an iterative search algorithm to find the solution subspaces. 


It is important to remark  here that an optimal solution can have less than $l$ subspaces, and the dimensions of the subspaces
can be less than $n$.
Since the minimization we consider is over  unions of no more than $l$ subspaces, where the dimension of the subspaces is no bigger than $n,$
some of the optimal solutions  for a given $(l,n)$ (that is, some of the solutions that give the smallest $\varepsilon$) 
will yield the minimum $l_0 \leq l$ such that
the data is $(l_0,n,\varepsilon)$-sparse, that is $l$ is set to be just an upper bound for the number of allowable subspaces.
Furthermore, the number $n$ constraining the dimension of the subspaces is also only an upper bound, that is, an optimal solution 
can have subspaces of dimension  strictly less than $n$. 



\subsection{ Organization and Contribution}
In this paper we solve the abstract  problem described in \eqref {total-error-0}. Unlike prior work in the subspace segmentation problem (see e.g., \cite {VMS05} and the references therein), it does not assume that the data $\F$ comes from union of subspaces, but instead it finds the best union of subspaces that matches the data, and therefore it is well adapted for subspace segmentation in the presence of noise, and for the problem of sparsity and dictionary  design in compressed sensing. Moreover, the setting includes finite and infinite dimensional spaces, and therefore can be used to solve the signal modeling problem described in \cite {LD07}. The subspaces that are sought are not restricted to be orthogonal, or with equal dimensions or with trivial intersection, and  there can be any number of subspaces up to a prescribed number $l$. 

In Section \ref {theory} we formally state as Problem 1 the question described by \eqref{total-error-0} and introduce a general abstract scheme for solving this specific problem, together with a lemma that will provide the tool for an algorithm described in a later section. 

In Section \ref {SISCONSM}  we consider the case of the Hilbert space $L^2(\R^d)$ and  where the infinite dimensional subspaces are  shift-invariant. We show that the general theory in Section \ref {theory} applies to this case and thereby we solve \eqref {total-error-0} in this situation. As a consequence we show how the signal modeling problem is solved in Section \ref {OSM}.

In Section \ref {R^N} we consider the finite dimensional case $\R^N$ and particularize the solution found in Section \ref {theory} to this case. This allows us to tie our method to the problem of finding sparsity models and the problem of finding optimal dictionaries as described in Section \ref {SOD}.

In Section \ref {ALG} we present an iterative search algorithm for finding the solution to Problem \eqref {total-error-0}. We prove that the algorithm terminates in finitely many steps. The algorithm is iterative and switches between subspace estimation and data segmentation in  a way that is similar to the subspace segmentation methods and the K-SVD method described in \cite {AEB06a, AEB06b,MYDF07,VMS05} and the references therein.


\section{Abstract Hilbert Space Case}\label{theory}

In this section we will introduce an abstract scheme that, in particular, contains the problems mentioned in the introduction. This scheme is much more general and can be used in many other situations. 

In this theoretical setting we will prove the existence of  optimal solutions and provide the
mathematical background for the algorithms to find these solutions.   

We will start by describing the basic ingredients for that setting and introducing some required notation.
First we will define the class of subspaces that we will use for the minimization.
 
Let $\HH$ be a Hilbert  space. For $x,y \in \HH$ let us denote by $d(x,y)=\|x-y\|_{\HH}$.
Given a finite subset $\F \subset \HH$ and a closed subspace $V$ of $\HH$, we denote by $E(\F,V)$ the total distance of the data set $\F$ to the subspace $V$, i.e.
\begin{equation}
\label {optimal0}
E(\F,V) = \sum_{f\in \F} d^2(f,V).
\end{equation}
We set $E(\F,V)=0$ for $\F=\emptyset$ and any  subspace $V$ of $ \HH$. 

Let   $\mathcal C$ be a  family of closed subspaces of $\HH$ containing the zero subspace.
We will say that $\CC$ has the Minimal Approximation Property  (MAP) if
for any finite set $\F$ of vectors in $\HH$ there exist a subspace $V_0 \in \CC$
that minimizes $
E(\F,V)$ over all the subspaces  $V \in \CC.$ That is,
\begin{equation}\label{optimal}
E(\F,V_0)=\min\limits_{V\in\CC}E(\F,V)\leq E(\F,V), \qquad \forall \;V \in \CC.
\end{equation}

Any subspace  $V_0 \in \CC$ satisfying \eqref{optimal}
 will be called an {\em optimal subspace for }$\F.$
Note that if $\F=\emptyset$ then every subspace in $\CC$ is optimal. We will choose
the zero subspace in that case.
For the rest of this section we will assume that the class $\CC$ has the Minimal Approximation Property.

Next, since we are interested in models that are union of subspaces, we will arrange the subspaces
 in finite bundles that will be our main objects, and define the distance (error) between a bundle and a set of vectors.

To do this, let us  fix $m, l \in \N$ with $1 \leq l \leq m$ and let $\F= \{f_1,...,f_m\}$ be a finite set of vectors in $\HH.$

\begin{figure}
\begin{center}
\includegraphics[width=\textwidth]{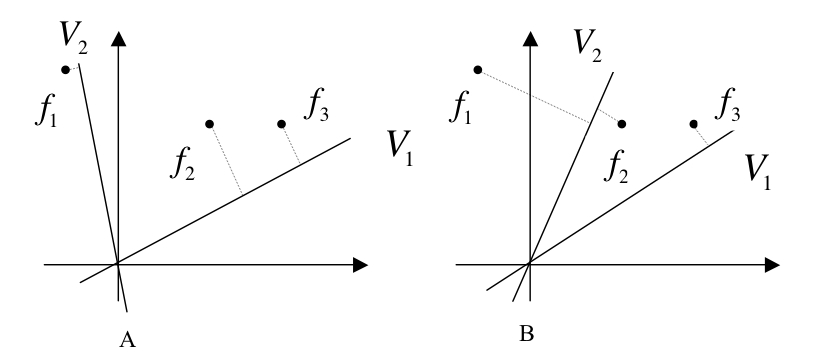}
\caption{\small Illustration for the objective function in  Problem~\ref {problem}:  A data set consists of three points $\F=\{f_1,f_2,f_3\}$ in  $\R^2$. A) Value of the objective function is $e=d(f_1,V_2)+d(f_2,V_1)+d(f_3,V_1)$; and   B)  Value of the objective function is $e=d(f_1,V_2)+d(f_2,V_2)+d(f_3,V_1)$. Note that the configuration of $V_1, V_2$  in Panel A forced a partition of the data into $P_1=\{f_1\}$ and $P_2=\{ f_2,f_3\}$, while the  configuration in B forced the partition  $P_1=\{f_1,f_2\}$ and $P_2=\{f_3\}$ for the same data. }
\label{FigNL}
\end{center}
\end{figure}

Define $\VV$ to be the set of sequences of elements in $\CC$ of length  $l$, i.e.
$$\VV=\VV(l)=\big\{ \EV = \{V_1,\dots,V_l\}:  V_i \in \CC, 1\le i\le l \big\}.$$
We will call these finite sequences {\it \word s}.
For $\EV \in \VV$  with $\EV = \{V_1,...,V_l\}$, we define,
\begin{equation}\label{total-error}
e(\F,\EV) = \sum_{f \in \F} \min_{1\leq j \leq l} d^2(f,V_j).
\end{equation} 
\begin {remark}
Note that $e(\F,\EV)$ is computed as follows: For each $f\in \F$ find the space $V_{j(f)}$ in $\EV$ closest to $f$,  compute $d^2(f, V_{j(f)})$, and then sum over all values  found by letting  $f$ run through $\F$ (see Figure \ref {FigNL}). Also note that $e(\F,\EV)$ is a non-linear function of $\F$.
\end {remark}
Hence, for the problems described in the introduction, what we want is to minimize $e$ over all possible \word s of subspaces.
This is formulated in the following problem.
\begin{prob}\label{problem}
\ 

\begin{enumerate}

\item Given a finite set $\F \subset \HH$, minimize $e(\F,\EV)$ over $\EV \in \VV$. 
That is, find
\beq
\inf \{ e(\F,\EV): \EV \in \VV\}.
\eeq
\item Find a \word\ $\EV_0 \in \VV$ (if it exists) such that  
\beq\label{solution} 
e(\F,\EV_0) = \inf \{ e(\F,\EV): \EV \in \VV\}.
\eeq 
\end{enumerate}
Any $\EV_0 \in \VV$ that satisfies \eqref{solution} will be called a {\em solution} to Problem \ref{problem}
\end{prob}
We will show (Theorem~\ref{main}) that Problem~\ref{problem} can be solved, i.e. for a given data set $\F$, there does exist a \word\ $\EV_0$ that minimizes $e(\F,\EV)$. Moreover, we propose an algorithm to find this optimal \word .

\begin{figure}
\begin{center}
\includegraphics[width=5cm]{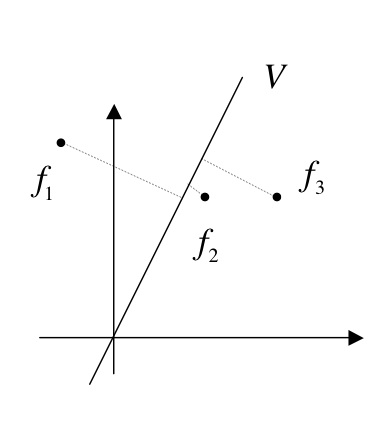}
\caption{\small Illustration for the objective function in Problem~\ref {problem}: Same data set $\F=\{f_1,f_2,f_3\}$ as in Figure \ref {FigNL}, but for a single subspace $V$. This objective function is the classical least squares cost function.}
\label{FigLin}
\end{center}
\end{figure}

If  we set  $l=1$, $\HH=\R^d$, and $\CC=\DD_n$ to be the set of all subspaces of dimension smaller (or equal) than $n$, then Problem~\ref {problem} reduces to the classical least squares problem.  This last  problem  has been studied extensively (see Figure \ref {FigLin} for an illustration in $\R^2$), and it can be solved using the well-known Singular Value Decomposition (SVD) (see e.g., \cite{Sch07, EY36}).

Before stating the main results, we need to give some definitions and set some notation.

We will denote by $\Pi =\Pi_l$  the set of all $l$-sequences $P= \{\F_1,\cdots, \F_l\}$ 
of subsets of $\F$ satisfying the property that
for all $ 1 \leq i,j \leq  l  $,
$$\F_i  \subset \F, \quad \F = \cup_{s=1}^{l} \F_s, \quad \text{ and }
\quad \F_i \cap \F_j =\emptyset \text{ for } i \neq j.$$

Note that we allow some of the elements of $P\in \Pi$ to be the empty set. By abuse of language
we will still call the elements of $\Pi_l$ {\em partitions} (of $\F$).  

For $P \in \Pi_l,\;  P = \{\F_1,...,\F_l\}$ and $\EV \in \VV, \EV=\{V_1,...,V_l\}$ we define,
\begin{equation}\label {Gamma}
\Gamma(P,\EV) = \sum_{i=1}^{l} E(\F_i,V_i).
\end{equation}
So $\Gamma$ measures the error between a fixed partition $P$ and a fixed \word\ $\EV$.

The following relations between partitions in $\Pi_l$ and \word s of length $l$ in $\CC$ will be relevant for our analysis.

Given a \word\  $\EV \in \VV, \;\; \EV =\{V_1,...,V_l\}\subset \CC$, we can partition the set $\F$ into a {\em best partition} $P = \{\F_1,...,\F_l\}$, by grouping together into $\F_i$, the vectors in $\F$ that are closer to a given subspace $V_i$ than to any other subspace $V_j$, $j\ne i$ (see Figure \ref {FigNL}). However, there are situations in which a vector $f\in \F$ is at equal distance from two or more subspaces from the bundle. Hence, there may be more than one partition associated to each $bundle$, and there is a subset $\Pil (\EV) \subset \Pi_l$ of {\em best partitions} in $\Pi_l$ naturally associated to $\EV$ defined by
\begin{quote}
$P=\{\F_1,...,\F_l\} \in \Pi_l$  is a member of $\Pil (\EV)$ if it satisfies\\
$f \in \F_j$ implies that  $d(f,V_j) \leq d(f,V_h) , h=1,...,l.$
\end{quote}

Conversely,  since $\CC$ has the MAP, given a partition $P=\{\F_1,...,\F_l\}$, we can define a  {\em best bundle} $ \EV_P=\{V_1,...,V_l\} \in \VV$ by finding (for each $i$) the space $V_i$ that minimizes \eqref {optimal0} for the given $\F_i$. However, there are situations in which, for some $\F_i$, there are more than one  subspace $V_i$ that minimizes  \eqref {optimal0}. Hence, there is a subset $\vv(P) \subset \VV$ of {\em best  \word s} associated to $\F$ defined by

\begin{quote}
$\EV_P = \{V_1,...,V_l\} \in \VV$ is a member of $\vv(P)$ if $V_i$ is an optimal subspace for $\F_i$
(in the sense of \eqref{optimal}) for each $i=1,...,l.$
\end{quote}

In what follows when we refer to a best partition associated to a \word\ $\EV$  we will mean, any element in $\Pil (\EV)$. Similarly, when we talk of a best \word\  associated to a partition $P$, this will mean an element in $\vv(P)$.

We also consider the set of all pairs $(P,\EV_P)$, where $P \in \Pi_l$ and $\EV_P \in \vv(P)$.
We will say that a pair $(P_0,\EV_{P_0})$ is {\em $\Gamma$-minimal} if
\begin{equation} \label{opt-part}
\Gamma(P_0,\EV_{P_0}) \leq \Gamma(P,\EV_P) \end{equation}
for all such pairs. 

Note that when trying to compute $e(\F,\EV)$, for each $f \in \F$ we first have to find the subspace $V_{j(f)}$ in $\EV$  that is closest to $f$ and then compute $d^2(f, V_{j(f)})$ (see remark after the definition of $e(\F,\EV)$ and Figure \ref {FigNL}). While for $\Gamma$, a partition is given and  we just compute the distance of each function to its corresponding space (not the closest one necessarily). The surprising fact is that $e$ and $\Gamma$ can indeed be compared, as the following lemma shows. In addition, this result
will later give us the key to obtain an algorithm for Problem \ref{problem}.

\begin{lemma}\label{leaf}
Let $(P_0, \EV_{P_0})$ be a $\Gamma$-minimal pair. Then we have 
\begin{equation}\label{equation-leaf}
e(\F,\EV_{P_0}) = \Gamma(P_0,\EV_{P_0}).
\end{equation}
\end{lemma}

\begin{proof}
It is clear that  $e(\F,\EV_{P_0}) \leq \Gamma(P_0,\EV_{P_0})$. 

For  the other inequality, if $\EV_{P_0} = \{V_1,...,V_l\}$ then for any $P \in \Pil (\EV_{P_0})$ we have
\begin{equation}
e(\F,\EV_{P_0}) = \sum_{i=1}^m \min_{1\leq j \leq l} d^2(f_i,V_j) = \Gamma(P,V_{P_0}).
\end{equation}

In addition, 
$\Gamma(P,V_{P_0}) \geq \Gamma(P,\EV_P),$ with $\EV_P \in \vv(P)$.
But by the minimality of  $\Gamma(P_0,V_{P_0})$ given by hypothesis, we have that
$ \Gamma(P,\EV_P) \geq \Gamma (P_0,\EV_{P_0})$, and the lemma follows.
\end{proof}

We are now ready to prove the following theorem which shows that we can solve Problem~\ref{problem}.

\begin{theorem}\label{main}
Let $\HH$ be a Hilbert space,  $m,l $ positive integers with $ l \le m$ and
$\F=\{f_1,...,f_m\}$ a set of vectors in $\HH$. Then 
\begin{enumerate}
\item  
There exists  a \word\ $\EV_0 \in \VV$
that solves Problem \ref{problem} for the data $\F$, that is,
$$ e(\F,\EV_0) = \inf \{ e(\F,\EV): \EV \in \VV\}.$$
\item 
If $(P_0, \EV_{P_0})$ is a $\Gamma$-minimal pair, then all the elements of $\vv(P_0),$ are solutions to Problem \ref{problem}.
\item Furthermore, if $\EV_0$ is a solution to Problem \ref{problem}, then there exists $P_0 \in \Pi_l$ such that $\EV_0 \in \vv(P_0)$, (i.e. $(P_0, \EV_0)$ is a $\Gamma$-minimal pair.
 \end{enumerate}
\end{theorem}
In other words, Theorem \ref {main} states that Problem~\ref {problem} has a solution for every finite set of vectors  $\F=\{f_1,...,f_m\}\subset \HH$ and every $l\ge 1$  if and only if $\CC$ has the MAP property.
One direction of the theorem is trivial. The interesting implication is that  if Problem~\ref {problem} can be solved for any $\F$ and $l=1$ then it can be solved for any $\F$ and any $l\ge1$. 

\begin{proof}

We will prove that if $(P_0, \EV_{P_0})$ is a $\Gamma$-minimal pair,   then 
$$e(\F,\EV_{P_0}) \leq e(\F,\EV), \quad \forall \; \EV \in \VV.$$

For this, let us choose an arbitrary $\EV \in \VV$. We have  that for each $P  \in \Pil (\EV)$
$$\Gamma(P,\EV) = e(\F,\EV).$$

Clearly $\Gamma(P,\EV_P) \leq \Gamma(P,\EV)$, for each $\EV_P \in \vv(P)$.

Because of the minimality of $\Gamma(P_0,\EV_{P_0})$, we have
$$\Gamma(P_0,\EV_{P_0}) \leq \Gamma(P,\EV_{P}).$$

As a consequence of  Lemma \ref{leaf} we know  that then
$$\Gamma(P_0,\EV_{P_0}) = e(\F,\EV_{P_0}),$$ which proves,
$$e(\F,\EV_{P_0}) \leq e(\F,\EV).$$
This shows that if $(P_0, \EV_{P_0})$ is a $\Gamma$-minimal pair, then each bundle $\EV_{P_0}$ solves Problem \ref{problem} for the data \F. Since the total number of pairs is finite, then there exist minimal pairs. This proves parts (1) and (2) of the Theorem.

For part (3) let $\EV_0 \in \VV$ be a solution to Problem ~\ref{problem}, i.e. $e(\F,\EV_0) \leq e(\F,\EV), \; \forall \; \EV \in \VV$. Consider $P_0 \in \Pil (\EV_0)$ and let $\EV_{P_0} \in \vv(P_0)$. Then, since $P_0 \in \Pil (\EV_0)$ and by the minimality of $\EV_0$ we have
$$ \Gamma(P_0,\EV_0) = e(\F,\EV_0) \leq e(\F,\EV_{P_0}) \leq \Gamma(P_0,\EV_{P_0}).$$
Therefore, $\Gamma(P_0,\EV_0) \leq \Gamma(P_0,\EV_{P_0})$, but by definition of $\Gamma$,
$\Gamma(P_0,\EV_{P_0}) \leq \Gamma(P_0,\EV)$ for any $\EV \in \VV$. So, 
$$\Gamma(P_0,\EV_0) = \Gamma(P_0,\EV_{P_0}), \quad \text{and} \quad \EV_0 \in \vv(P_0).$$
Moreover, $(P_0, \EV_0)$ is $\Gamma$-minimal since
$$ \Gamma(P_0,\EV_0) = e(\F,\EV_0) \leq e(\F,\EV_{P}) \leq \Gamma(P,\EV_{P}).$$
This completes the proof of the Theorem.
\end{proof} 

\begin{remark}
If $0 < l_1 < l_2$, then for any $\EV \in \VV(l_1)$, $\EV = \{V_1, \dots, V_{l_1}\}$ the \word\ $\EV' = \{V_1, \dots, V_{l_1}, \{0\}, \dots, \{0\}\}$ belongs to  $\VV(l_2)$ and therefore, we have
$$e(\F,\EV_{P_0}(l_1)) \geq e(\F,\EV_{P_0}(l_2)).$$
So the error decreases (or  at least does not increase) when $l$ (the number of subspaces) increases. Note that in case that the number of subspaces equals the number of data, the error is zero, since we can pick for each data signal the subspace spanned by itself.

It is important to remark here that optimal \word s can have the zero subspace as some of its components.
So, if $l_0$ is the number of subspaces that have dimension greater than zero, in some optimal \word \; $\EV_0$,  then
the \word \, with $l_0$ components obtained after the $l-l_0$ zero components are removed from $\EV_0,$ is also an optimal 
\word \, for the Problem \ref{problem} when $\VV(l)$ is replaced by $\VV(l_0)$. Thus as mentioned in the introduction, the number $l$ is simply a set to be an a priori upper bound on the number of subspaces, and the optimal solution(s) can have any number of subspaces $l_0\le l$.

\end{remark}


\section{The Shift-Invariant Space Case and Optimal Nonlinear Signal Models}
\label {SISCONSM}
In this section we will apply the theory of Section \ref{theory} to the Hilbert space  $L^2(\R^d)$. In order to do that we will select a family of subspaces with the Minimal Approximation Property. We will  describe in what follows the necessary setting.

We begin by recalling the definition of  frames and some of their properties  (see for example \cite {Cas00, Chr03, Gro01, HW96}). 

 Let $\HH$ be a Hilbert space and $\{u_i\}_{i\in I}$ a countable subset of $\HH$. The set $\{u_i\}_{i\in I}$ is said to form a {\em frame}  for $\HH$ if there exist $q,Q>0$ such that 
 \[
 q\|f\|^2\le \sum\limits_{i\in I}|<f,u_i>|^2\le Q \|f\|^2, \quad \forall \, f \in \HH.
 \]
 If $q=Q$, then $\{u_i\}_{i\in I}$ is called a {\em tight frame}, and it is called a {\em Parseval frame} if $q=Q=1$.

 If $\{u_i\}_{i\in I}$ is a Parseval frame for a subspace $W$ of a Hilbert space $\HH$, and if $a \in \HH$, then the orthogonal projection of  $a$ onto $W$ is given by:
  \begin{equation} \label{proj}
\proj_W(a) = \sum_{i\in I} \langle a, u_i \rangle u_i.
\end{equation}
Thus, a Parseval frames acts as if it were an orthonormal basis of $W$, even though it may not be one.

\subsection{Shift-Invariant Spaces}
\label {SIScase}
 In this paper, a shift-invariant space will be a  subspace of $L^2(\R^d)$ of  the form:
\begin{equation}
\label {sis}
S(\Phi):=\hbox{closure}_{L_2} \;  \hbox{span} \{ \varphi_i(x-k): \;  i =1,\ldots, n, \, k \in \Z^d\} 
\end{equation}
       where $\Phi=\{\varphi_1,..., \varphi_n\}$ is a  set of functions in
$L^2(\R^d)$.        The  functions $\varphi_1 ,
\varphi_2, \ldots, \varphi_n$  are called a  {\em set of generators} for the
space $S=S(\Phi)=S(\varphi_1,\ldots,\varphi_n)$  and any such space $S$ is called a {\em finitely generated shift-invariant space (FSIS)} (see e.g.,
\cite {Bow00}). These spaces are often used as standard signal and image models. For example, if
$n=1,\ d=1$ and
$\phi(x) = \hbox {sinc}(x)$, then the underlying space is the
space of \index{band-limited functions} band-limited functions (often
used in communications).

Finitely generated shift-invariant spaces,
 can have different sets of generators. The {\em length} of an FSIS  $S$ is,
$$
l(S)=\min  \{\ell \in \N: \exists \; \varphi_1,\ldots,\varphi_\ell \in S \text{ with }  S=S(\varphi_1,\ldots,\varphi_\ell)\}.
$$
If $S = \{0\},$ we set $l(S) =0.$
We will denote by $\mathcal L_n$  the set of all  shift-invariant spaces with length less than or equal to $n$. That is, an element in $\mathcal L_n$ is a \SIS \  that has a set of $s$ generators with $s \leq n$.

\subsection{The Minimal Approximation Property for SIS}

In \cite{ACHM07} it was proven that $\LL_n$ has the MAP. More precisely, 
\begin {theorem}\label{best-fit}
\label {existence}
Let $\F=\{f_1,\dots,f_m\}$ be a set of functions in $L^2(\R^d)$. Then

 there exists $V\in \LL_n$ such that 
\begin{equation}\label{best}
\sum \limits_{i=1}^m \|f_i-\proj_Vf_i\|^2\le \sum \limits_{i=1}^m \|f_i-\proj_{V'}f_i\|^2,  \quad \forall \ V' \in \LL_n.
\end{equation}

\end {theorem}
Here $\proj_V$ denote the orthogonal projection onto the subspace $V.$

Furthermore, an explicit description of an optimal space (that is not necessarily  unique) 
and an estimation of the error, was obtained in \cite{ACHM07}, as is described below.
Let us call,
\begin {equation}
 \mathcal E (\F,n)=\min\limits_{V'\in \mathcal L_n}\sum \limits_{i=1}^m \|f_i-\proj_{V'}f_i\|^2.
\end {equation} 
To compute the error $ \mathcal E (\F,n)$ we need to consider the Gramian matrix $G_{\F}$ of 
$\F=\{f_1,\dots,f_m\}.$ Specifically, the {\em Gramian} $G_\Phi$ of a set of functions $\Phi=\{\varphi_1,\dots,\varphi_n\}$ with elements in $L^2(\R^d)$ is defined to be the $n \times n$ matrix of $\Z^d$-periodic functions 

\begin{equation} \label{gram}
[G_{\Phi}(\zo)]_{i,j} = \sum_{k\in\Z^d}\widehat \varphi_i(\zo+k)\overline{\widehat \varphi_j(\zo+k)}, \qquad \zo\in \R^d,
\end{equation}
where $\widehat \varphi_i$ denotes the Fourier transform of $\varphi_i$, and $\overline{\widehat \varphi_i}$ denotes the complex conjugate  of $\widehat \varphi_i$.


The next  theorem produces a set of generators for an optimal space $V \in \mathcal L_n$ and   provides a formula for the exact value of  the error. 
 \begin {theorem} \label{teo-gram} \cite{ACHM07} Under the same assumptions as in Theorem \ref {existence}, let   $\lambda_1(\zo)\ge\lambda_2(\zo)\ge\dots\ge \lambda_m(\zo)$ be the eigenvalues of the Gramian $G_{\F}(\zo)$.  Then
\label {construction}
\begin {enumerate}
\item The eigenvalues $\lambda_i(\zo) $, $1\le i \le m$ are $\Z^d$-periodic, measurable functions in $L^2([0,1]^d)$ and 
\begin{equation} 
 \mathcal E(\F,n)= \sum \limits_{i=n+1}^m\int\limits_{[0,1]^d}\lambda_i(\zo)d\zo .
\end{equation}

\item \label{teo-4} Let $E_i:=\{\zo\;:\; \lambda_i(\zo)\ne 0\}$, and define $\tilde \sigma_i(\zo)=\lambda_i^{-1/2}(\zo)$ on $E_i$ and $\tilde \sigma_i(\zo)=0 $ on $E_i^c$. Then, there exists a choice of measurable left eigenvectors   $y_1(\zo),\dots,y_n(\zo)$ with $y_i=(y_{i1},...,y_{im})^t, i=1,...,n,$ associated with the first $n$ largest eigenvalues of $G_{\F}(\zo)$ such that the functions defined by
\begin{equation}\label{generators}
\hat \varphi_i(\zo)=\tilde \sigma_i(\zo)\sum \limits_{j=1}^my_{ij}(\zo)\hat f_j(\zo), \quad i=1,\dots,n, \; \zo\in \R^d
\end{equation}
are  in $L^2(\R^d)$.

 Furthermore, the corresponding set of  functions $\Phi=\{\varphi_1,\dots,\varphi_n\}$ is a set of generators for an optimal space $V$ and the set $\{\varphi_i(\cdot-k), k \in \Z^d, i=1,\dots,n\}$ is a Parseval frame for $V$. 
\end{enumerate}
 \end {theorem}
 
Note that \eqref{generators} says that in particular the generators of the optimal space are $l_2(\Z)$-linear combinations of the integer translates of the data $\F$.

\subsection{Best Approximation by Bundles of SIS}
\label {BASIS}

Let  $\F= \{f_1,...,f_m\}$  be functions in $L^2(\R^d)$  and   $n$ a positive integer smaller than $m$.

The result in Theorem \ref{best-fit} says that the class $\LL_n$ has the Minimal Approximation Property.

Define $\SiS=\{\{S_1,...,S_l\}: S_i \in \LL_n\}$ to be the set of \word s
of SIS in $\LL_n.$
Now we can apply Theorem \ref{main} to conclude that

\begin{theorem}
\label {OPSIS}
Let $\F=\{f_1,...,f_m\}$ vectors in $L^2(\R^d)$, then there exist a \word\
$\ES_0 = \{S_1^0,...,S_l^0\} \in \mathcal S$ such that
\begin{equation}
e(\F,\ES_0) = \sum_{i=1}^m  \min_{1\leq j \leq l} d_2^2(f_i,S_j^0)
\leq  \sum_{i=1}^m  \min_{1\leq j \leq l} d_2^2(f_i,S_j)
\end{equation}
over all \word s $\ES=\{S_1,\dots,S_l\} \in \SiS.$

\end{theorem} 

Let  $P_0=\{\F^0_1,...,\F^0_l\}$ be a best partition of $\F$ associated to the optimal \word\  
$\ES_0=\{S_1^0,\dots,S_l^0\}$
(i.e. $P_0=\{\F^0_1,\dots,\F^0_l\} \in \Pil (\ES_0)$, is such that $f_j \in \F^0_i$ implies $d(f_j,S_i^0) \leq d(f_j,S_k^0), \quad k=1,...,l).$ 

Using Theorem \ref{teo-gram} for each $h =1,...,l$ and such that $\F^0_h\neq \emptyset,$ a set of generators forming a Parseval frame can be obtained for the optimal space $S_h^0$ in terms of the singular values and singular vectors of the gramian $G_{\F_h}$ associated to the subset $\F^0_h.$
Furthermore a formula for the minimum error $E(\F^0_h,S_h^0)$ is given in terms of the singular values of  $G_{\F^0_h}$. Therefore, thanks to Lemma \ref {leaf} and the definition of $\Gamma$ (see \eqref {Gamma}),    $e(\F,\ES_0)$ can be computed exactly.

\subsection {Optimal Signal Models}
\label {OSM}
If it is known a priori that a class of signals belong to a union of shift-invariant spaces  $\cup^l_{i=1}S_i$ each of length no larger than $n$, then combining Theorems \ref {teo-gram} and   \ref  {OPSIS}, we can find the signal model $\cup^l_{i=1}S_i$ exactly and the generators of each space $S_i$. This solution includes the case where the signal class consists of a single shift-invariant space solved in \cite {ACHMR03}. 

When the data $\F$ is corrupted by noise, or if the true signals are not from a union of shift-invariant spaces  $\cup^l_{i=1}S_i$ of length no larger than $n$, but we still wish to model the class of signals by such a union, then we can use the \word\ $S_0=\{S^0_1,\cdots,S^0_l\}$  found in Theorem \ref  {OPSIS} to obtain an optimal  signal model  $\cup_i S^0_i$  compatible with the observed data. The optimal signal model $\cup_iS^0_i$ is a union of  infinite dimensional spaces $S_i^0$, but each $S_i^0$ is a shift-invariant space that can be generated by  at most $n$ frame generators  $\Phi_0^i=\{\varphi^i_{0,1},\dots,\varphi^i_{0,s_i}\}$, $s_i\le n$, $i=1,\dots, q$, $q\le l$ (we only use the spaces $S^0_i$ that have length larger than $0$). Each signal $f_j\in \F$ can now be modeled by its orthogonal projection $f_j^a=\proj_{S^0_i}f_j$ onto its closest space $S^0_i$. which consists of countably many but generally infinite linear combinations of atoms, i.e., $f^a_j=\sum\limits^{s_i}_{p=1}\sum_k c^p_k\varphi^i_{0,p}( .  -k)$.

Note that if the data $\F$ is corrupted by noise, then the optimal model can be used as a denoising method. Other applications of the optimal signal model are those of learning, data segmentation, and classification.


\section{The finite dimensional case  $\R^N$, sparse representations, and optimal dictionaries}\label{R^N}

In this section we will consider  Problem \ref{problem} for the case in which the Hilbert space is $\R^N$ (in applications, usually one thinks of  $N$ as being very large). In this case, our data $\F= \{f_1,\dots,f_m\}$ are vectors in $\R^N$. 
Let us denote by $\DD_n$ the set of all subspaces of dimension smaller (or equal) than $n$. 
To see that $\DD_n$ has the MAP, we will recourse to some well-known results about Singular Value Decomposition (SVD) and in particular, the Eckart-Young Theorem. 

Let  us briefly recall the SVD decomposition, (for a  detailed treatment see for example \cite{HJ85}, or the Appendix of \cite{ACHM07}).
Let $A \in \R^{N\times m}$ with columns $ \{a_1, \dots, a_m\}$, and let $r$ be the rank of $A$. One can obtain its SVD as follows.  Consider the matrix
  $A^tA \in \R^{m\times m}$. Since $A^tA$ is self-adjoint and positive semi-definite, 
 its eigenvalues $\lambda_1\ge \lambda_2\ge \cdots \ge \lambda_m$   are nonnegative and the associated eigenvectors $y_1,\ldots, y_m$ can be chosen to form an orthonormal basis of $\R^m$.   Note that the rank $r$ of $A$ corresponds to the largest index $i$ such that $\lambda_i>0$. 
  The left singular vectors $u_1,\ldots, u_r$ can then be obtained from 
$$\sqrt{\lambda_i}u_i=Ay_i, \;\;  \text{ that is } u_i = \lambda_i^{-1/2} \sum _{j = 1}^m y_{ij} a_j.
\qquad (1\le i\le r).$$
Here $y_i=(y_{i1},...,y_{im})^t.$
The remaining left singular vectors $u_{r+1},\ldots,  u_m$ can be chosen to be any orthonormal collection  of
$m-r$ vectors in $\R^N$ that are perpendicular to $\text{ span } \{a_1,\ldots, a_m\} $.   One may then readily verify that  
 \begin{equation}
 \label{SVD}
 A=\sum_{k=1}^m \sqrt{\lambda_k} u_k y_k^t = U\Lambda^{1/2} Y^t ,
  \end{equation}
  where $U \in \R^{N\times m}$ is the matrix $U = \{u_1, \dots, u_m\},
 \;\Lambda^{1/2} = \text{diag} (\lambda_1^{1/2},...,\lambda_m^{1/2}),$
 and  $Y = \{y_1,...,y_m\} \in \R^{m \times m}$
  with $U^tU=I_m = Y^tY=YY^t.$
 
The following theorem of Schmidt (cf. \cite{Sch07}) (usually coined as Eckart-Young Theorem\cite{EY36}) shows that our set $\DD_n$ has the MAP. (We again denote by $\proj_V$ the orthogonal projection onto the space $V$.) 

\begin {theorem}[Eckart-Young] \label{EY}
Let $\{f_1, \dots, f_m\}$ be a set of vectors in $\R^N$ and  $r=$ dim (span$\{f_1, \dots, f_m\}). $ 
Suppose that the associated matrix $A =\{f_1,\ldots, f_m\}$, has SVD $A=U\Lambda^{1/2} Y^t$ and that $0<n\le r$.  
 If  $W=\sspan \{u_1,\dots,u_n\}$, then
$$\{\proj_{W}f_1,\ldots, \proj_{W}f_m\} =\sum_{i=1}^n \sqrt{\lambda_i}u_iy_i^t   =A_n$$ and
\begin {equation}
\sum_{i=1}^m \|a_i-\proj_{W}a_i\|_{2}^2 \le \sum_{i=1}^m \|a_i-\proj_{V}a_i\|_{2}^2, \quad \forall \; V \in \DD_n.
\end {equation}
Furthermore,  the   space $W$ is unique if $\lambda_{n+1}\neq \lambda_n$. In addition, 
\begin{equation}\label{err-ey}
{E}(\F,W) = \min\limits_{V\in \DD_n}\sum \limits_{i=1}^m \|f_i-\proj_{V}f_i\|^2 = \sum_{j=n+1}^r \lambda_j . \end{equation}
\end {theorem}

\subsection{Best Non-Linear Approximation by Bundles of Subspaces in $\R^N$}
\label{map-rn}

Let $\F = \{f_1, \dots, f_m\}$ be a set of vectors in $\R^N$ and $n \leq m$. As indicated before, Theorem~\ref{EY}  states precisely that  $\DD_n$ has the MAP.

Define again $\BB = \BB(l)$ to be the set of non-empty \word s of length $l$ in $\DD_n$, i.e. 
$$
\BB = \{\{V_1,...,V_l\}: V_i \in \DD_n,\, \;\; i=1,\dots,l\}.
$$
 We then have the following theorem.

\begin{theorem}\label {OPRN}
Let $\F = \{f_1, \dots, f_m\}$ be vectors in $\R^N$, and let $l$ and $n$ be given ($l < m$, $n < N$), then there exist a \word\ $\EV_0 = \{V_1^0, \dots, V^0_l\} \in \BB$, such that
$$ e(\F,\EV_0) = \sum_{i=1}^n \min_{1\leq j\leq l} d_2^2(f_i, V^0_j) = \inf \{ e(\F,\EV): \EV \in \BB\}.
$$ 
\end{theorem}

Let   $P_0=\{\F^0_1,...,\F^0_l\}$ be the best partition of $\F$ associated to the optimal \word\  
$\EV_0 = \{V_1^0,...,V_l^0\}$
(i.e. $P_0 = \{\F^0_1,\cdots,\F^0_l\} \in \Pil (\EV_0)$, is such that $f_j \in \F^0_i$ implies $d(f_j,V_i^0) \leq d(f_j,V_k^0), \; k=1,...,l).$ 

Now, using Theorem \ref{EY} for each $h =1,...,l$ and such that $\F^0_h\neq \emptyset,$ a set of generators forming an orthonormal base can be obtained for the optimal space $V_h^0$ in terms of the singular values and singular vectors of the matrix $A_h$ associated to the subset $\F^0_h.$
Furthermore, by \eqref{err-ey}, a formula for the minimum error $E(\F^0_h,V_h^0)$, and therefore for  $e(\F,\EV_0)$, is given in terms of the singular values of  $A_h$.

\begin {remark}
Theorem \ref  {OPRN} remains true if we replace the set $\DD_n$ by the
set $\DD_{(n_1,\dots,n_l)}$ of bundles $\{V_1,\dots,V_l\}$ such that $\dim
V_i\le n_i$, for $i=1,\dots,l$.
\end {remark}

\subsection {Sparsity and Optimal Dictionaries}
\label {SOD}

Now, we will describe the relation between the solution to Problem \ref{problem} for $\R^N,$ as described in Section \ref{map-rn}, with the problem of dictionary finding and sparsity.
Let us introduce the following problem. See for example \cite{AEB06b} . 
\begin{prob}\label{problem-2}

 Given data $\F=\{f_1,...,f_m\}$ in $\R^N$ and positive integers $n$ and $d$, find a dictionary $D$ (i.e. a set of vectors in $\R^N$) of length at most $d,$ such that each  $f_i$ can be written as a linear combination of at most $n$ atoms in $D.$
 
 That is, (in matrix notation) find a $N\times r$ matrix  $D$, with columns $a_1,...,a_r$ in $\R^N$ and $r \leq d$,  such that there exist an $r \times m$ matrix $X$, with columns $x_1,...,x_m$ in $\R^r$,
  such that $\F = DX$ with $\|x_i\|_0 \leq n,$ for $ i=1,...,m.$
\end{prob}

Now we will introduce a definition of sparsity and  show its connection with Problem \ref{problem-2}.
\begin{definition} \label{def-sparse}
Let $n,l, m$ be positive integers, with $n, l < m.$

Given a set of vectors $\F=\{f_1,...,f_m\}$ in $\R^N$
and a real number $\varepsilon \geq 0$, we will say that the data $\F$ is $(l,n,\varepsilon)$-{\it sparse}
if there exist subspaces $V_1,...,V_l,$ of $\R^N$ with dim$(V_i)\leq n$ for $i=1,...,l$, such that 
 \begin{equation}
e(\F,\big \{V_1,...,V_l \big \}) = \sum_{i=1}^m \min_{1\leq j \leq l} d^2( f_i,V_j) \leq \varepsilon.
\end{equation}
When $\F$ is $(l,n,0)$-sparse, we will simply say that $\F$ is $(l,n)$-sparse.
We will also say that the data is $\varepsilon$-sparse if the values of $l$ and $n$ are clear from the context.
\end {definition}
Note that if $\F$ is $(l,n,\varepsilon)$-sparse, then it is also $(l,n,\eta)$-sparse for every $\eta \geq \varepsilon.$ So, usually it is interesting to know the minimun $\varepsilon$ such that the data is
$\varepsilon$-sparse.
The above definition of sparsity is an intrinsic property of the data and the Hilbert space in which the data lives, and does not depend on any specific dictionary. 

\subsubsection{The case $\varepsilon = 0$}
Let us now consider the case $\varepsilon = 0$.

If the data $\F$ is $(l,n)$-sparse and  $V^0_1,...,V^0_l$ are optimal spaces 
(that is when
$ e(\F,\left \{V^0_1, \dots, V^0_l \right \})= 0$ and dim$(V^0_i)\leq n$ for $i=1,...,l$)  then each $f \in \F$ belongs to some of the spaces $\{V^0_i\}_{i=1,...,l}.$

 For each $i=1,...,l,$ let us call $r_i = dim(V^0_i)$ 
and let $\{w_{i1},...,w_{ir_{i}}\}$ be an orthonormal basis of $V^0_i$.
Define 
$$D= \bigcup_{i=1}^l \{w_{ij}:  1 \leq j \leq r_{i} \} \subset \R^N.$$

The vectors in $D$ have the property that each $f \in \F$  can be written as a linear combination of at most $n$ elements in $D$.
In other words, we have found a dictionary $D$  such that it solves Problem \ref{problem-2},
for the data $\F$ and length $s$ with $s = r_1+\dots+r_l \leq ln$.

So Theorem \ref{OPRN} provides a solution  to Problem \ref{problem-2} with a dictionary of length at most $ln$, in case that the data is $(l,n)$-sparse. We will see below that if the data $\F$ is not $(l,n)$-sparse, then Theorem \ref{OPRN},
provides the minimum $\varepsilon$ such that the data is $(l,n,\varepsilon)$-sparse.

Note that if the  basis of each $V^0_i$ is properly chosen  then in many cases the number of atoms in the dictionary can be reduced, due to the fact that the subspaces can have non-trivial intersections.

So, as before, let $V^0_1,...,V^0_l$ 
be optimal spaces, and let $\mathcal{U}=\{u_1,...,u_s\}$ be a set of vectors with the property that for each $i \in \{1,...,l\}$ there 
is a subset  $\mathcal{U}_i \subset \mathcal{U}$ such that span$(\mathcal{U}_i)= V^0_i$. 
Set $D_0 = \{w_1,...,w_{s_0}\}$ to  be a minimal set with this property.

Then Theorem \ref{OPRN} implies that $D_0$ is a dictionary  that solves Problem
\ref{problem-2} for data $\F$ and positive integers $n$ and $d=ln$.
We want to remark here that a minimal set is not a linearly independent set in general.
It is not difficult to see that considering all possible intersections of the subspaces $V^0_i$, a minimal set can be constructed.



\subsubsection{The case $\varepsilon > 0$}

If the data is not $(l,n)$-sparse, then Theorem \ref {OPRN} implies that there is no dictionary $D$ of length $d\ =  ln$ or smaller  that solves Problem \ref{problem-2},
 for the data $\F$.

If we still want to find a dictionary  of length no larger than $ln$ with $\|x_i\|_0 \leq n$
for $i = 1,...,n$, then the question is: what error do we have to allow in order to have a solution? In other words, 
what is the minimum $\varepsilon$ such that the data $\F$ is $(l,n,\varepsilon)$-sparse?
This question gives rise to  the following extension of Problem \ref{problem-2}.

\begin{prob}
\label {problem-3}
Let $\F=\{f_1,...,f_m\}$ in $\R^N$, and  $n$, $d$ positive integers.
With the same notation as in Problem \ref{problem-2}, given $\varepsilon \geq 0$, find a dictionary $D$ with no more than $d$ atoms and a matrix $X$ such that
$$\|\F-DX\| \leq \varepsilon$$ with $\|x_i\|_0 \leq n,$ for $ i=1,...,m.$
\end{prob}
Theorem \ref{OPRN} provides in this case the  solution for  the minimum possible error and establishes the exact value of the error. More precisely,
let $V^0_1,...,V^0_l$ be a \word\ of  optimal subspaces and let $\varepsilon=e(\F,\big \{V^0_1,...,V^0_l \big \})$. Let us choose as before a minimal set $D_0=\{w_1,...,w_{s_0}\}$ for that solution,
then we have that there exist vectors $x_1,...,x_m$ in  $R^{s_0}$ such that
$$\|\F-D_0X\| = \varepsilon $$ with $\|x_i\|_0 \leq n,$ for $ i=1,...,m.$
Furthermore, given any $N \times r$ matrix $D$ with $r \leq s_0$ and any
matrix $X$ with columns $x_1,...,x_m$ in $\R^r$ and $\|x_i\|_0 \leq n$ for $ i=1,...,m$, we have
$$ \|\F-DX\| \geq \varepsilon.$$ 

So, a solution of Problem \ref{problem} gives an optimal solution for Problem \ref{problem-3}
and finds the exact (optimal) $\varepsilon$-sparsity. 

Note that when $n$ or $l$ increase then in general the minimum error will decrease.
It is also important to emphasize here that some subspaces from an optimal \word \, for the data $\mathcal{F}$ and $(l,n)$
can have dimension zero, so in applications these subspaces can be removed and the non-trivial subspaces will produce the same error. Thus, the subspaces that are found are not restricted to be orthogonal, or with equal dimensions or with trivial intersection, and  there can be any number of subspaces up to a prescribed number $l$, and each subspace can be of any dimension up to  a prescribed number $n$.


\section {Search Algorithm}
\label {ALG}
Although  Theorem \ref {main} establishes the existence of a global minimizer solution to Problem \ref {problem}, an exhaustive search over all possible partitions is not feasible in practice and a search algorithm is needed. Lemma \ref{leaf} used in the proof  Theorem \ref {main} suggests an iterative search algorithm that we will present in this section, and we will show that the algorithm always terminates in finitely many steps. The search algorithm for finding the solution to Problem \ref {problem} is given in the program below, with the notation of Section \ref{theory}. It uses two choice functions $\CW, \CP$, where $\CW$ is a choice function assigning $ P \longmapsto \EV_P\in \vv(P)$, and $\CP$ is a choice function assigning $ \EV \longmapsto P\in \Pil (\EV)$.

\begin{algo} \
\label {Alg}
\begin {enumerate}
\item Pick any partition $P_1\in \Pi_l$;
\item Find and choose $\EV_{P_1}=\CW(P_1) \in \vv(P_1) \subset \VV$ by minimizing $\Gamma(P_{1},\EV)$ over $\EV\in \VV$;
\item Set $j=1$;
\item {\bf While} $\Gamma(P_j,\EV_{P_j})> e(\F,\EV_{P_j})$ ;
\item Choose  a new partition $P_{j+1}=\CP(\EV_{P_j})\in \Pil (\EV_{P_j})$ associated to $\EV_{P_j}$;   
\item Find and choose $\EV_{P_{j+1}}=\CW(P_j) \in \vv(P_j)$, by minimizing $\Gamma(P_{j+1},\EV)$ over $\EV\in \VV$;
\item Increase $j$ by $1$, i.e., $j\rightarrow j+1$;
\item {\bf End while}
\end {enumerate}
\end{algo}
 Note that this algorithm, starting from a bundle $\EV_{P_1}$ in step (2), produces a sequence of bundles $\EV_{P_1}, \EV_{P_2}, \EV_{P_3}, \dots$ with the property that  that $e(\F,\EV_{P_1})\geq e(\F,\EV_{P_{2}}) \geq e(\F,\EV_{P_{3}}) \geq \dots $. The algorithm stops precisely when  for some $j\geq 1$ $e(\F,\EV_{P_j}) =e(\F,\EV_{P_{j+1}})$. We will now  see that the algorithm terminates in finitely many steps.
\begin{proof}
 We first note that if $\Gamma(P_j,\EV_{P_j})> e(\F,\EV_{P_j})$, then $P_{j+1}\notin \Pil (V_{P_j})$. To see this, we argue by contradiction: if $P_{j+1}\in \Pil (V_{P_j})$, then $\Gamma(P_{j+1},\EV_{P_{j+1}})= \Gamma(P_j,\EV_{P_j})$. But we  have that 
 $$\Gamma(P_{j+1},\EV_{P_{j+1}})\le e(\F,\EV_{P_j})< \Gamma(P_j,\EV_{P_j}),$$
  which is a contradiction.
 
  Therefore, since the set of partitions is finite, the algorithm must stop in at most $\#\Pi_l$ steps. The algorithm terminates when $\Gamma(P_{end},\EV_{P_{end}})=  e(\F,\EV_{P_{end}})$.
\end{proof}
\begin{note}
A  partition $P_{m}$ such that $\Gamma(P_{m},\EV_{P_{m}})=  e(\F,\EV_{P_m})$ will be called a {\it  minimal  partition} (see remark below).
\end{note}
 
\begin {remark}
\end {remark}
\begin {enumerate}
\item The algorithm can be formulated as a search for minimal partitions in the partially ordered set $(\Pi_l, \preceq)$, where the order of the elements in $\Pi_l$ depends on the specific choice functions $\CW,\CP$ in (2), (5) and (6) of the algorithm.
Specifically, $P \preceq Q$ if there exists an integer $s\ge 1$ such that $P =(\CP\CW)^sQ$. Since  $(\Pi_l, \preceq)$ is a partially ordered set with finitely many elements, a nonempty set of minimal  partitions $\mathcal M \subset \Pi_l$ exists.

\item The algorithm  will always terminate in finitely many steps but the bundle $\EV_{P_m}$ associated to  the final minimal  partition $P_m$ may not be the global minimizer of Problem \ref {problem}.  The algorithm can be viewed as a search in a directed graph whose vertices are the partitions. Each iteration moves from one partition to the next via a directed edge. If the graph has a single component, then the algorithm will always end at a  partition whose associated bundle is a global minimizer. However, if the graph has more than one component, then the algorithm will end up at a  partition, whose associated bundle is not necessarily the global minimizer. 



\item In the search algorithm, steps (2) and (6) must be implemented by some other minimizing algorithms. For the two cases that we studied in this paper, a space $\EV_{P_{j+1}}$ that minimizes $\Gamma(P_{j+1},\EV)$ over $\EV \in \VV$ can be explicitly found and computed, by the Eckhard-Young Theorem for subspaces of $\HH=\R^N$, and by  Theorem 2.1 (\cite{ ACHM07}) for shift-invariant spaces of $\HH=L^2(\R^d)$. 
 Both methods are based on the Singular Value Decomposition, they are easily implemented, and all the approximation errors can be computed  exactly. 


\item In searching for a dictionary with $d$ atoms such that each data point $f\in \R^N$ from a set of data $\mathcal{F}=\{f_1,\dots,f_m\}$ is approximated by a single atom, our method coincides with the K-SVD algorithm proposed in (\cite{AEB06a, AEB06b}) and produces the same dictionary if steps (2) and (6) are implemented using the SVD. However, for the case where each data point is approximated by a linear combination of $n>1$ atoms, the two methods are not comparable, even if steps (2) and (6) are implemented using the SVD.
\end {enumerate}

\section{Conclusions}

Theorem \ref {main} can be viewed as a way of finding an optimal  (generaly non-linear) signal model of the form $\cup^l_{i=1}S_i$ from some observed data. For example, the application of Theorem \ref {main} to shift-invariant spaces in Section \ref {SIScase} produced Theorem \ref  {OPSIS},  which gives the optimal  signal model of at most $l$ shift-invariant spaces compatible with the observed data. The resulting solution is an optimal \word\ $S_0=\{S^0_1,\cdots,S^0_l\}$ that consists of a finite sequence of  infinite dimensional spaces, such that each space $S_i^0$ of this sequence is generated by  the integer translates of at most $n$ generators. This type of best signal model $\cup_{i\in I}S^0_i$ derived from a set of observed data $\mathcal{F} = \{f_1, \dots, f_m\} \subset L^2$ may be used for example in sampling and reconstruction as well as other applications.

If applied to $\R^N$, Theorem \ref {main}, can be viewed as a way of finding an optimal  sparse representation of the data $\F \subset \R^N$, optimal dictionaries, and subspace segmentation in the presence of noise as discussed in Section \ref {SOD}. Specifically, the application of Theorem \ref {main} to the finite dimensional case $\R^N$ produces Theorem \ref {OPRN} which gives the optimal $\varepsilon$-sparse representation of the data as discussed in Section \ref {SOD}.  In particular, Theorem \ref{OPRN} proves the existence of a dictionary with minimal error and minimal length for sparse data representation.

One contribution of this work is that it unifies and complements some of the new non-linear techniques used in sampling theory, the Generalized Principle 
Components Analysis, and the dictionary design problem. However, there are still many questions that need to be addressed before this methodology becomes applicable. For example,  the algorithm proposed in the last section may end up at a local minimum which is not a global minimum. The termination of the algorithm depends on the initial condition. Thus, an interesting question is to estimate the number of minima in terms of some characteristics 
of $\CC$ and $\F$. Another interesting question is to estimate the speed at which the algorithm converges in 
terms of some characteristics of $\CC$ and $\F$. Testing  the dependence of algorithm on the initial partition, the data, and noise level, for the case $\HH=\R^N$ and  $\CC=\DD_n$ using an SVD implementation for steps (2) and (6)   is also important for future research, and applications.






\begin{thebibliography}{BDDW07}

\bibitem[ACHM07]{ACHM07}
A.~Aldroubi, C.~A. Cabrelli, D.~Hardin, and U.~M. Molter, \emph{Optimal shift
  invariant spaces and their parseval frame generators}, Applied and
  Computational Harmonic Analysis (2007).

\bibitem[ACHMR03]{ACHMR03}
Akram Aldroubi, C.~Cabrelli, D.~P. Hardin, U.~Molter, and E.~Rodado,
  \emph{Determining sets of shift invariant spaces}, Proceedings of ICWA
  (Chenai, India), 2003.

\bibitem[AEB06a]{AEB06b}
M.~Aharon, M.~Elad, and A.M. Bruckstein, \emph{The k-svd: An algorithm for
  designing of overcomplete dictionaries for sparse representation}, IEEE
  Trans. On Signal Processing \textbf{54} (2006), no.~11, 4311 -- 4322.

\bibitem[AEB06b]{AEB06a}
Michal Aharon, Michael Elad, and Alfred~M. Bruckstein, \emph{On the uniqueness
  of overcomplete dictionaries, and a practical way to retrieve them}, Linear
  Algebra Appl. \textbf{416} (2006), no.~1, 48--67. \MR{MR2232919
  (2007a:94026)}

\bibitem[AG01]{AG01}
A.~Aldroubi and K-H. Gr{\"o}chenig, \emph{Non-uniform sampling in
  shift-invariant space}, Siam Review \textbf{43} (2001), no.~4, 585--620.

\bibitem[BDDW07]{BDDW07}
R.~Baraniuk, M.~Davenport, Ronald~A. DeVore, and M.~Wakin, \emph{A simple proof
  of the restricted isometry property for random matrices}, Preprint, 2007.

\bibitem[Bow00]{Bow00}
Marcin Bownik, \emph{The structure of shift-invariant subspaces of
  ${L}^2(\mathbb{R}^n)$}, Journal of Functional Analysis \textbf{177} (2000),
  282--309.

\bibitem[Cas00]{Cas00}
Peter~G. Casazza, \emph{The art of frame theory}, Taiwanese J. Math. \textbf{4}
  (2000), no.~2, 129--201.

\bibitem[Chr03]{Chr03}
Ole Christensen, \emph{An introduction to frames and {R}iesz basis}, Applied
  and Numerical Harmonic Analysis, Birkh\"auser, 2003.

\bibitem[CR06]{CR06}
E.~Cand{\`e}s and J.~Romberg, \emph{Quantitative robust uncertainty principles
  and optimally sparse decompositions}, Foundations of Comput. Math. \textbf{6}
  (2006), 227--254.

\bibitem[CRT06]{CRT06}
E.~Cand{\`e}s, J.~Romberg, and Terence Tao, \emph{Robust uncertainty
  principles: Exact signal reconstruction from highly incomplete frequency
  information}, IEEE Trans. on Information Theory \textbf{52} (2006), 489--509.

\bibitem[CT06]{CT06}
E.~Cand{\`e}s and Terence Tao, \emph{Near optimal signal recovery from random
  projections: Universal encoding strategies}, IEEE Trans. on Information
  Theory \textbf{52} (2006), 5406--5425.

\bibitem[DeV07]{Dev07}
Ronald~A. DeVore, \emph{Deterministic constructions of compressed sensing
  matrices}, Preprint, 2007.

\bibitem[Don06]{Don06}
David~L. Donoho, \emph{Compressed sensing}, IEEE Trans. on Information Theory
  \textbf{52} (2006), 1289--1306.

\bibitem[DVB07]{DVB07}
Pier~Luigi Dragotti, M.~Vetterli, and T.~Blu, \emph{Sampling moments and
  reconstructing signals of finite rate of innovation: Shannon meets
  strang-fix}, IEEE Transactions on Signal Processing \textbf{55} (2007),
  1741--1757.

\bibitem[EY36]{EY36}
C.~Eckart and G.~Young, \emph{The approximation of one matrix by another of
  lower rank}, Psychometrica \textbf{1} (1936), 211 -- 218.

\bibitem[GN03]{GN03}
R.~Gribonval and M.~Nielsen, \emph{Sparse decompositions in unions of bases},
  IEEE Trans. Inf. Theory \textbf{49} (2003), 3320--3325.

\bibitem[Gr{\"o}01]{Gro01}
K.~Gr{\"o}chenig, \emph{Foundations of time-frequency analysis}, Applied and
  Numerical Harmonic Analysis, Birkh{\"a}user, 2001.

\bibitem[HJ85]{HJ85}
R.~Horn and C.~Johnson, \emph{Matrix analysis}, Cambridge University Press,
  Cambridge, 1985.

\bibitem[HW96]{HW96}
Eugenio Hern{\'a}ndez and Guido Weiss, \emph{A first course on wavelets}, CRC
  Press, Boca Raton, FL, 1996.

\bibitem[LD07]{LD07}
Y.~Lu and M.~N. Do, \emph{A theory for sampling signals from a union of
  subspaces}, IEEE Transactions on Signal Processing, (2007).

\bibitem[MYDF07]{MYDF07}
Y.~Ma, A.~Yang, H.~Derksen, and R.~Fossum, \emph{%
Estimation of Subspace Arrangements with Applications in Modeling and Segmenting Mixed Data},
 preprint, to appear in SIAM Review, 2007.

\bibitem[MDHW07]{MDHW07}
Y.~Ma, H.~Derksen, W.~Hong and J.~Wright, \emph{%
 Segmentation of Multivariate Mixed Data via Lossy Coding and Compression},
IEEE Transactions on Pattern Analysis and Machine Intelligence (PAMI), \textbf{29}, no. 9, (2007), 1546--1562.

\bibitem[MV05]{MV05}
I.~Maravic and M.~Vetterli, \emph{Sampling and reconstruction of signals with
  finite rate of innovation in the presence of noise}, IEEE Transactions on
  Signal Processing \textbf{53} (2005), 2788--2805.

\bibitem[RSV06]{RSV06}
Holger Rauhut, K.~Schass, and P.~Vandergheynst, \emph{Compressed sensing and
  redundant dictionaries}, Preprint, 2006.

\bibitem[Sch07]{Sch07}
E.~Schmidt, \emph{Zur theorie der linearen und nichtlinearen
  integralgleichungen. i teil. entwicklung willk{\"u}rlichen funktionen nach
  system vorgeschriebener}, Math. Ann. \textbf{63} (1907), 433--476.

\bibitem[Tro04]{Tro04}
J.~A. Tropp, \emph{Greed is good: Algorithmic results for sparse
  approximation}, IEEE Trans. Inf. Theory \textbf{50} (2004), 2231--2242.

\bibitem[VMS05]{VMS05}
R.~Vidal, Y.~Ma, and S.~Sastry, \emph{Generalized principal component analysis
  (gpca)}, IEEE Transactions on Pattern Analysis and Machine Intelligence
  \textbf{27} (2005), 1--15.

\end{thebibliography}
\newcommand{\etalchar}[1]{$^{#1}$}
\providecommand{\bysame}{\leavevmode\hbox to3em{\hrulefill}\thinspace}
\providecommand{\MR}{\relax\ifhmode\unskip\space\fi MR }
\providecommand{\MRhref}[2]{%
  \href{http://www.ams.org/mathscinet-getitem?mr=#1}{#2}
}
\providecommand{\href}[2]{#2}

\end{document}